\documentclass[11pt,leqno,twoside, headinclude]{amsart}
\usepackage{enumerate}

\title{On fibrations with formal elliptic fibers}
\def\titl{On fibrations with formal elliptic fibers}
\def\auth{Manuel Amann and Vitali Kapovitch}
\date{December 13th, 2011}
\thanks{The first named author was supported by a Research Grant of the German Research Foundation.
The second author was supported in part by a Discovery Grant from NSERC}
\author{\auth}

\usepackage[colorlinks,pdftex, plainpages=false]{hyperref}
\hypersetup{pdftitle=\titl, pdfauthor=\auth, pdftoolbar=false,
plainpages=false, hyperindex=true, pdfdisplaydoctitle=true}

\hypersetup{colorlinks,%
citecolor=black,%
filecolor=black,%
linkcolor=black,%
urlcolor=black,%
pdftex}

\usepackage{color}

\usepackage{amsmath, amssymb, amscd, amsthm}
\usepackage{stmaryrd}
\usepackage{latexsym}
\usepackage[all]{xy}
\usepackage{pb-diagram}
\usepackage{rotating}
\usepackage{multicol}
\usepackage{lscape}
\usepackage{wasysym}
\usepackage{longtable}
\xyoption{all}

\newtheorem{theo}{Theorem}[section]
\newtheorem{main}{Theorem}
\newtheorem{maincor}[main]{Corollary}
\newtheorem*{main*}{Theorem}

\newtheorem{mainconj}{Conjecture}
\newenvironment{mainconjec}[1]{\begin{mainconj}[#1]\normalfont}{\end{mainconj}}

\newtheorem{defi2}[theo]{Definition}
\newenvironment{defi}{\begin{defi2}\normalfont}{\end{defi2}}
\newenvironment{defin}[1]{\begin{defi2}[#1]\normalfont}{\end{defi2}}
\newtheorem{rem2}[theo]{Remark}
\newenvironment{rem}{\begin{rem2}\normalfont}{\hfill$\boxbox$\end{rem2}}
\newtheorem{lemma}[theo]{Lemma}
\newtheorem{cor}[theo]{Corollary}
\newtheorem*{cor*}{Corollary}

\newtheorem*{conj*}{Conjecture}
\newtheorem*{theo*}{Theorem}
\newtheorem*{ques*}{Question}
\newtheorem*{mi2}{Main Idea}

\newtheorem{ex2}[theo]{Example}
\newenvironment{ex}{\begin{ex2}\normalfont}{\hfill$\boxbox$\end{ex2}}
\newtheorem{exer2}[theo]{Exercise}

\newtheorem{alg2}[theo]{Algorithm}

\newcommand{\cc}{{\mathbb{C}}}                                     % complex numbers
\newcommand{\kk}{{\mathbb{K}}}                                     % arbitrary field
\newcommand{\hh}{{\mathbb{H}}}                                     % quaternions
\newcommand{\qq}{{\mathbb{Q}}}                                     % rational numbers
\newcommand{\rr}{{\mathbb{R}}}                                     % real numbers
\newcommand{\pp}{{\mathbf{P}}}                                     % polynomials, projective space
\newcommand{\Gr}{{\mathbf{Gr}}}                                    % Grassmannian
\newcommand{\s}{{\mathbb{S}}}                                      % sphere
\newcommand{\zz}{{\mathbb{Z}}}                                     % integers
\newcommand{\SO}{{\mathbf{SO}}}                                    % special orthogonal group
\newcommand{\U}{{\mathbf{U}}}                                      % unitary group
\newcommand{\SU}{{\mathbf{SU}}}                                    % special unitary group
\newcommand{\Sp}{{\mathbf{Sp}}}                                    % symplectic group
\newcommand{\B}{{\mathbf{B}}}                                      % Lie group type B
\newcommand{\E}{{\mathbf{E}}}                                      % Lie group type E
\newcommand{\F}{{\mathbf{F}}}                                      % Lie group type F
\newcommand{\G}{{\mathbf{G}}}                                      % Lie group type G
\newcommand{\Spin}{{\mathbf{Spin}}}                                % spin group
\newcommand{\dif} {{\operatorname{d}}}                             % differential operator d
\newcommand{\In} {{\,\subseteq\,}}                                 % subset
\newcommand{\Ni} {\,{\supseteq}\,}                                 % superset
\newcommand{\Aut}{{\operatorname{Aut}}}                            % automorphisms
\newcommand{\End}{{\operatorname{End}}}                            % endomorphisms
\newcommand{\Hol}{{\operatorname{Hol}}}                            % Holonomy group
\newcommand{\id}{{\operatorname{id}}}                              % identity
\newcommand{\APL}{{\operatorname{A_{PL}}}}                         % polynomial differential forms
\newcommand{\cat} {{\operatorname{cat}}}                           % Lusternik--Schnirelman category
\newcommand{\cupl}{{\operatorname{c}}}                             % Cup-length
\newcommand{\cl}{{\operatorname{cl}}}                              % Cone-length
\newcommand{\e}{{\operatorname{e}}}                                % Toomer invariant
\newcommand{\rk}{{\operatorname{rk\,}}}                            % toral rank
\newcommand{\diag}{{\operatorname{diag}}}                          % diagonal matrix
\newcommand{\Der}{{\operatorname{Der}}}                             % imaginary part

\newcommand{\comment}[1]{}                                         % insert a large comment
\newcommand{\xto}[1]{\xrightarrow{#1}}                             % abbreviation for \xrightarrow
\newcommand{\hto}[1]{\overset{#1}{\hookrightarrow}}                % abbreviation for labelled \hookrightarrow
\newcommand{\rto}[1]{\overset{#1}{\rightarrow}}                % abbreviation for labelled \rightarrow

\newcommand{\biq}[2]{#1\;\!\!\!\sslash \;\!\!\!#2}                 % creates a biquotient
\newcommand{\case}[1]{\textbf{Case #1.}}                           % Case
\newcommand{\ack}{\noindent\textbf{Acknowledgements. }}            % Acknowledgements
\newcommand{\str}{\noindent\textbf{Structure of the article. }}    % Structure of the article
\def\co{\colon\thinspace}      % colon in maps
\def \dD{\mathcal D}    %differential on derivations
\newenvironment{prf}{\begin{proof}[\textsc{Proof}]} {\end{proof}}     % alternative environment proof
\begin{document}

\maketitle \thispagestyle{empty}

%%%%%%%%%%%%%%%%%%%%%%%%%%%%%%%%%% Abstract%%%%%%%%%%%%%%%%%%%%%%%%%%%%%%%%%%%%%%%%

\begin{abstract}
We prove that for a fibration of simply-connected spaces of finite type $F\hto{} E\to B$ with $F$ being positively elliptic and $H^*(F,\qq)$ not possessing non-trivial derivations of negative degree, the base $B$ is formal if and only if the total space $E$ is formal. Moreover, in this case the fibration map is a formal map. As a geometric application we show that positive quaternion K\"ahler manifolds are formal and so are their associated twistor fibration maps.
\end{abstract}

%%%%%%%%%%%%%%%%%%%%%%%%%%%%%%%%%% Introduction %%%%%%%%%%%%%%%%%%%%%%%%%%%%%%%%%%%

\section*{Introduction}

The problem we shall address in this article gains its appeal  from two rather disjoint sources; the first one being  inherent to algebraic topology, the other one motivated by a prominent question in Riemannian geometry.
%We shall approach the geometric issue via giving a special answer to the topological question within the realm of Rational Homotopy Theory.

Let us begin by illustrating the motivating  geometric setting.
Riemannian manifolds with special holonomy form a very interesting class of spaces which include
 \emph{K\"ahler manifolds} (manifolds with $\U(n)$-holonomy), \emph{Calabi--Yau
manifolds} (manifolds with $\SU(n)$-holonomy) and \emph{Joyce manifolds} \linebreak[4]($\Spin(7)$-holonomy and $\G_2$-holonomy).
We will be interested in \emph{quaternion K\"ahler manifolds}, which are manifolds
with  holonomy
contained in $\Sp(n)\Sp(1)$. Such manifolds are known to be Einstein and they are called \emph{positive}, if
their scalar curvature is positive.

Positive quaternion K\"ahler geometry lies in the intersection of very classical yet
rather different fields in mathematics. Despite its geometric
setting, %which involves fundamental definitions from Riemannian geometry,
 it was  discovered to be accessible by methods from
(differential) topology, symplectic geometry and complex algebraic
geometry.

To the knowledge of the authors, the approach by rational homotopy theory, which we provide in this article, is the first
one of its kind in the setting of quaternion K\"ahler geometry.

The field of positive quaternion K\"ahler geometry settles around the following
\begin{mainconjec}{LeBrun, Salamon}\label{conj01}
Every positive quaternion K\"ahler manifold is a symmetric space.
\end{mainconjec}
%Several prominent articles give partial answers to this classification problem
There are a number of partial results supporting this conjecture
and only symmetric examples, the so-called \emph{Wolf spaces} are known. However, the conjecture remains open in general. Thus our motivating geometric question will be

\begin{ques*}[Geometry]
How close are positive quaternion K\"ahler manifolds to being symmetric?
\end{ques*}

\vspace{5mm}

Let us now describe the topological motivation, which arises from rational homotopy theory. This  is a very elegant and easily-computable
version of homotopy theory at the expense of losing information on
torsion. It provides a transition from topology to algebra by
encoding the rational homotopy type of a space in a commutative
differential graded algebra. In particular, rational homotopy groups
as well as Massey products can be derived from the algebra
structure. Likewise, the rational cohomology algebra of the space is
the homology algebra of the corresponding commutative differential
graded algebra.

The concept of formality features prominently amongst the properties
of topological spaces, as this property reduces the study of the
rational homotopy type entirely to the problem of merely
understanding the rational cohomology algebra. Or in other words: We
may derive the rational cohomology from the rational homotopy type,
however, is the information contained in the rational cohomology
already sufficient to reconstruct the rational homotopy type? If the
answer is ``yes", the space is called \emph{formal}.

Although it is not known to the authors if  the following conjecture is stated explicitly in the literature, it is widely believed that
\begin{mainconjec}{}
A simply-connected compact Riemannian manifold of special holonomy is a formal space.
\end{mainconjec}
A famous result by Berger---having undergone several refinements---states that a simply-connected irreducible non-symmetric Riemannian manifold $M$ has one of the holonomy groups $\SO(n)$---the generic case not comprised in the term ``special holonomy'', with $\dim M=n$---or $\U(n)$ ($\dim M=2n$), $\SU(n)$ ($\dim M=2n$), $\Sp(n)$ ($\dim M=4n$), $\Sp(n)\Sp(1)$ ($\dim M=4n$), $\G_2$ ($\dim M=7$) respectively $\Spin(7)$ ($\dim M=8$).

We draw the attention of the reader to how nicely the prerequisites of irreducibility and being non-symmetric fit the context of formality: The finite Cartesian product of simply-connected spaces is formal if and only if so is each factor. Due to a famous result by Cartan, symmetric spaces are formal.

A celebrated result in \cite{DGMS75} states that compact K\"ahler manifolds are formal. Obviously, this comprises manifolds whose holonomy group is one of $\U(n)$, $\SU(n)$ and $\Sp(n)$. However, no further results in proving formality in the setting of special holonomy---yet a lot of attempts in that direction (!)---are known to the authors.

It is noteworthy that the formality of K\"ahler manifolds is indeed a geometric result and cannot be attributed to known topological properties of K\"ahler manifolds. In particular, there are compact simply-connected manifolds having the Hard--Lefschetz property but lacking formality (for example see \cite[Example 4.4 p.~346]{Cav07}). The result in \cite{DGMS75} is derived from the famous $\dif\dif^c$-lemma, which itself reverberates strongly in many generalisations in the literature, be it in the symplectic context, the generalized complex case etc.~(see \cite{Cav05})---always with the intent to relate it to formality. However, although all sorts of partial results and counterexamples exist, the $\dif\dif^c$-lemma often seems to be related to a Lefschetz-like structure.

As for our purposes let us just mention the attempt in \cite{Ver11} to generalize the $\dif\dif^c$-lemma to the context of special holonomy (and $\G_2$-manifolds, in particular). (Note that Joyce manifolds also satisfy a Lefschetz-like property.) The author's goal was to present a sufficient criterion to prove formality in this context. However, there seems to be a problem with the main tool, Proposition 2.19 on p.~1008 of that article. We are indebted to Spiro Karigiannis for pointing this out to us. Lastly, let us  mention that there exists an example of a simply-connected compact manifold sharing all the known topological properties of a $\G_2$-manifold but lacking formality---see \cite[Example 8.5, p.~131]{Cav05}.

The main geometric application of the topological  results of this article is to settle the conjecture above for positive quaternion K\"ahler manifolds, thus stating a first positive result in this direction, since the appearance of the article \cite{DGMS75}.

\vspace{5mm}

The following natural question---widely discussed in the literature---will be our topological motivation:

\begin{ques*}[Topology]
How do the formality of the base space and the one of the total space relate in a fibration?
\end{ques*}

The answer we shall provide to this question relies on the following concepts:
Recall that a fibration  $F\hto{} E\to B$ is called \emph{totally non-cohomologous to zero} or TNCZ for short, if the induced
homomorphism $H^*(E,\qq)\to H^*(F,\qq)$ is surjective. This is easily seen to be equivalent to  the Leray--Serre spectral sequence of this fibration degenerating on the $E_2$-term. It is also well-known that $H^*(F,\qq)$ has no negative degree derivations if and only if any fibration over a simply-connected base with fiber $F$ is TNCZ. We remark that an example of such fibres is provided by any simply-connected space whose cohomology algebra  satisfies  hard-Lefschetz duality. This applies, in particular, to all K\"ahler manifolds and (up to a degree shift, i.e.~using the Kraines form of degree $4$ in place of the K\"ahler form) to---per se rationally $3$-connected, but eventually to all---positive quaternion K\"ahler manifolds (cf.~\cite{Mei83}, \cite{Bla56}).

Moreover, recall that a simply-connected topological space $F$ is called \emph{positively elliptic} or \emph{$F_0$}, if it is rationally elliptic, i.e.~it has finite dimensional rational homotopy and cohomology, and if it has positive Euler characteristic. In this case its rational cohomology is concentrated in even degrees only. These spaces admit pure Sullivan models and feature prominently in rational homotopy theory. Classical examples of $F_0$-spaces are \emph{biquotients} (and, in particular, \emph{homogeneous spaces}) $\biq{G}{H}$ with $\rk G=\rk H$.

For the convenience of the reader we shall briefly review the notion of a biquotient. Let $G$ be a compact connected Lie group and let
$H\In G\times G$ be a closed Lie subgroup. Then $H$ acts on $G$ on
the left by $(h_1,h_2)\cdot g=h_1gh_2^{-1}$. The orbit space of this
action is called the \emph{biquotient} $\biq{G}{H}$ of $G$ by $H$.
 If the action of $H$ on $G$ is free, then
$\biq{G}{H}$ possesses a manifold structure. This is the only case
we shall consider. If $H=K\times L$ where $K\subset G\times 1$ and $L\subset 1\times G$ then the biquotiuent $\biq{G}{(K\times L)}$ is often denoted by $K\backslash G/L$.

Clearly, the category of biquotients contains the
one of homogeneous spaces. It was shown in \cite{Kap} that biquotients admit pure models \linebreak[4](cf.~\cite[Theorem 3.50]{FOT08}).

The most prominent conjecture which deals with $F_0$-spaces is
\begin{mainconjec}{Halperin}\label{conj02}
Suppose $F$ is an $F_0$-space.
Then  $H^*(F,\qq)$ has no negative degree derivations.
\end{mainconjec}
The conjecture holds true for  large classes of positively elliptic spaces:
It is satisfied for  homogeneous spaces ~\cite{ST87}, if the cohomology algebra $H^*(F,\qq)$ has at most $3$ generators
~\cite{Tho81, Lup90}, if all the generators are of the
same degree ~\cite{SZ75, SZ75a}, in the ``generic case''
~\cite{PP96} or---as already mentioned---in the case of hard-Lefschetz spaces.

Because of this conjecture we shall refer to spaces whose rational cohomology algebras do not possess negative degree derivations as spaces satisfying the Halperin conjecture---even if the spaces in question are \emph{not} rationally elliptic or of positive Euler characteristic.

It is a known fact (see \cite[Theorem 3.4]{Lup98}) that if $F$ is $F_0$ and satisfies Halperin's conjecture, then, given a fibration $F\hto{} E \to B$, the formality of the base space implies the formality of the total space.  In fact, more generally,  if a fibration is TNCZ and the fiber is  formal and elliptic, then the formality of the base implies the formality of the total space \cite[Proposition 3.2]{Lup98}.

Recall that a space $X$ is called of \emph{finite type} if all its cohomology groups over $\qq$ are finite dimensional and all its rationalised homotopy groups are finite dimensional.
Our main result and our proposed answer to the topological question is

%\newpage
\begin{main}\label{thm A'}
Let
\begin{align*}
F\hto{} E \rto{f} B
\end{align*}
be a fibration of simply-connected topological spaces of finite type. Suppose that   $F$ is elliptic,  formal and   satisfies the Halperin conjecture. Then  $E$ is formal if and only if $B$ is formal.

Moreover, if $B$ and $E$ are formal, then the map $f$ is formal.
\end{main}

It is well-known that $F_0$ spaces are formal; in fact, they are even \emph{hyperformal} and therefore \emph{intrinsically formal} (see  \cite[Remark 2.7.11.2, p.~120]{AP93}). Moreover, they admit pure Sullivan models.

Consequently, the following corollary is a direct consequence of Theorem \ref{thm A'} and completes the picture for $F_0$-spaces. However, we remark that Example \ref{ex01} given below shows that the class of formal, elliptic spaces satisfying Halperin's conjecture is strictly larger than the class of $F_0$-spaces satisfying the Halperin conjecture.
\begin{maincor}\label{theoA}
Let
\begin{align*}
F\hto{} E \rto{f} B
\end{align*}
be a fibration of simply-connected topological spaces of finite type. Suppose that $F$ is an $F_0$ -space which satisfies the Halperin conjecture. Then  $E$ is formal if and only if $B$ is formal.

Moreover, if $B$ and $E$ are formal, then the map $f$ is formal.
\end{maincor}

%\comment{%
%
%\begin{main}\label{theoA}
%Let
%\begin{align*}
%F\hto{} E \to B
%\end{align*}
%be a fibration of simply-connected topological spaces. Suppose $E$ and $F$ are  formal and $F$ satisfies Halperin's conjecture.
%
%Then $B$ is also formal.
%
%\end{main}
%}

%\comment{
%Combining these results with Theorem~\ref{theoA} we immediately obtain
%}
%As explained above Theorem \ref{theoA} applies in particular to all $F_0$-spaces satisfying the Halperin conjecture which conjecturally includes all %$F_0$-spaces.
%We do not know  whether either direction of the theorem remains true for non-elliptic fibers satisfying Halperin's conjecture.
We remark that there exists an example due to Thomas (see ~\cite[Example III.13]{Tho82}) of a TNCZ fibration (even a cohomologically trivial one!) with formal base and fiber but non-formal total space. In Example \ref{ex-lup} (due to Lupton) we produce a TNCZ fibration of simply-connected spaces with formal yet hyperbolic fiber satisfying Halperin's conjecture, formal base space and non-formal total space.

%Thus we suspect that the ``if" part of the above theorem is likely false in that generality.

\vspace{5mm}

The interlink between topology and geometry in our case is provided by the \emph{twistor fibration}
\begin{align*}
\cc\pp^1 \hto{} Z \to M
\end{align*}
of a Positive Quaternion K\"ahler Manifold $M$. (It is known that a positive quaternion K\"ahler manifold $M$ is compact and simply-connected.)

As the name ``quaternion K\"ahler''  indicates the subject lies in certain proximity to the field of K\"ahler
geometry. Indeed, (positive) quaternion K\"ahler manifolds can be
considered a quaternionic analogue of (compact) K\"ahler
manifolds. The twistor fibration is one way of illustrating this proximity, as the \emph{twistor space} $Z$ of a positive quaternion K\"ahler manifold is a K\"ahler manifold $Z$. K\"ahler
manifolds have been found to be formal spaces by joint work of
Deligne, Griffiths, Morgan and Sullivan~\cite{DGMS75}. Obviously, $\cc\pp^1\cong \s^2$ satisfies the Halperin conjecture.

Compact symmetric spaces are known to be formal; thus the formality of positive quaternion K\"ahler manifolds would be a consequence of a confirmation of Conjecture \ref{conj01}. Thus we can offer one more piece of the puzzle described via the geometric question we posed:
\begin{main}\label{maincor}
A positive quaternion K\"ahler manifold is formal. So is its twistor fibration.
\end{main}

Let us note that the only geometric input we use for this result is that positive quaternion K\"ahler manifolds are compact and simply-connected and that the total space of the twistor bundle admits a K\"ahler metric and hence is formal. Since this only uses Theorem~\ref{thm A'} with $F=\mathbb S^2$ one may suspect that in that special case it is easy to deduce formality of the base from formality of the total space. However, we do not  know of  a substantially simpler proof of Theorem~\ref{thm A'} even for $F=\mathbb S^2$ than the one we present for a general $F$.

 One may hope that some further progress in approaching the conjecture of LeBrun and Salamon (Conjecture \ref{conj01}) that positive quaternion K\"ahler manifolds are symmetric spaces can be obtained using methods from rational homotopy theory in conjunction with some further geometric input.

For example, one can try to improve Corollary \ref{maincor} in the following two directions, which might also be considered a motivation for proving the formality of positive quaternion K\"ahler manifolds in the first place:
%\comment{
%Let us  suggest and sketch two interpretations that might
%serve as a motivation for proving formality in the case of positive quaternion K\"ahler manifolds.

 On the one hand
formality is an obstruction to \emph{geometric formality}
(see \cite{Kot01}) which means that the product of
harmonic forms is harmonic again. Geometric formality enforces
strong restrictions on the topological structure of the underlying
manifold. For example, the Betti numbers of the manifold $M^n$ are
restricted from above by the Betti numbers of the $n$-dimensional
torus~\cite[Theorem 6]{Kot01}. (This result was  further strengthened
for K\"ahler manifolds by Nagy~\cite[Corollary 4.1]{Nag06}).
Symmetric spaces are geometrically formal (~\cite{KT03},
\cite{St02}). Thus it is tempting to conjecture the same for
Positive Quaternion K\"ahler Manifolds.

On the other hand the Bott conjecture speculates that
simply-connected compact Riemannian manifolds with nonnegative
sectional curvature are rationally elliptic. In the quaternionic
setting there are a number of results (see \cite[Theorem A, p.~150]{Fan04},
 \cite[Formula 14.42b, p.~406]{Bes08}) suggesting that positive
scalar curvature might be regarded as a substitute for positive
sectional curvature to a certain extent.

In the case of positive quaternion K\"ahler manifolds---mainly
since the rational cohomology is concentrated in even degrees
only---we suggest to see formality as a very weak substitute for
ellipticity. Indeed, if positive quaternion K\"ahler manifolds were
elliptic spaces---e.g.~like simply-connected homogeneous
spaces---then they would be $F_0$-spaces, which are formal.

If one
is willing to engage with this point of view, the formality of
positive quaternion K\"ahler manifolds may be seen as heading
towards a quaternionic Bott conjecture.

\vspace{5mm}

Let us briefly state the following trivial consequences of Corollary \ref{maincor}.
\begin{maincor}%\label{RHTcor03}
A rationally $3$-connected positive quaternion K\"ahler manifold $M^{4n}$
satisfies
\begin{align*}
n=\frac{\dim M}{4}=\cupl_0(M)=\e_0(M)=\cat_0(M)=\cl_0(M)
\end{align*}
If $M$ is not rationally $3$-connected, then $M\cong\Gr_2(\cc^{n+2})$ by Theorem~\ref{PQKtheo01} below and hence
\begin{align*}
2n=\frac{\dim M}{2}=\cupl_0(M)=\e_0(M)=\cat_0(M)=\cl_0(M)
\end{align*}
\end{maincor}
For the definition of the numerical invariants involved see the
definition on \cite[28, p.~370]{FHT01}, which itself relies on various
definitions on the pages 351, 360 and 366 in \cite[27]{FHT01}.

In order to prove this corollary one uses Theorems \ref{PQKtheo01} and \ref{PQKtheo02} cited below in order to reduce the problem to a rationally $3$-connected manifold. A volume form is given by $[u^n]$, where $u$ is the Kraines form in degree $4$. Due to rational $3$-connectedness, this allows to compute the rational cup-length $\cupl_0(M)=n$. The rest of the equalities then is a consequence of formality (cf.~\cite{FHT01}, Example 29.4, p.~388).

Recall that if $b_2(M)\neq 0$, we obtain that $M$ is a complex Grassmanian. It is K\"ahlerian, in particular, with K\"ahler form $\omega$.  The equation $\cupl_0(M)=2n$ follows from the existence of the volume form $\omega^{2n}$.

%\vspace{5mm}

 A
simply-connected quaternion K\"ahler manifold with vanishing scalar
curvature is hyper-K\"ahlerian and K\"ahlerian, in particular. So it
is a formal space. The twistor fibration in this case is the canonical projection $M\times \s^2\to M$. This directly yields
\begin{maincor}%\label{cor04}
A compact simply-connected non-negative quaternion K\"ahler manifold
is formal. The twistor fibration is formal.
\end{maincor}

Let us end the introduction with some more remarks:
In general, positive quaternion K\"ahler manifolds are not
\emph{coformal}, i.e.~their rational homotopy type is not
necessarily determined by their rational homotopy Lie algebra or,
equivalently, their minimal Sullivan models do not necessarily have
 strictly quadratic differentials. An obvious counterexample is
$\hh\pp^n$ for $n\ge2$.

In low dimensions, i.e.~in dimensions $12$ to $20$, relatively
simple proofs of  formality of  positive quaternion K\"ahler manifolds can be given using the concept of
$s$-formality developed in \cite{FM05}. Alternatively, one may use
the existence of isometric $\s^1$-actions on $12$-dimensional and
$16$-dimensional positive quaternion K\"ahler manifolds and further
structure theory to apply Corollary \cite[Theorem 5.9, p.~2785]{Lil03} which
yields formality.

\vspace{5mm}

\str In section \ref{sec01} we shall give a very brief introduction
to positive quaternion K\"ahler geometry whilst we do the same for the necessary techniques from rational homotopy theory in section \ref{sec02}.
%Sections
%\ref{sec03} and \ref{sec04} feature the proofs of theorems
%\ref{theoC} and \ref{theoD} respectively. In section \ref{sec05} the
%main results, theorems \ref{theoA} and \ref{theoB} will be derived.
Section \ref{sec03} is devoted to the proof of the Main Theorem \ref{thm A'}.
Finally, in section \ref{sec04}, we conclude with a depiction of
several counterexamples for possible statements similar to
Theorem \ref{thm A'} when assumptions on the fiber are weakened.
%---in the case of fibrations which are not TNCZ.

\vspace{5mm}

\textbf{As a general convention for this article we shall assume all spaces involved to be simply-connected (and, in particular, connected) and have finite type. Also all graded algebras we consider are assumed to be connected and also to have finite type, i.e.~they have finite dimensional cohomology in every dimension and are finitely generated in every degree.  Besides, cohomology is taken with rational coefficients and all commutative differential graded algebras are algebras over $\qq$.}

\vspace{5mm}

%\str In section \ref{sec01}

%\vspace{3mm}

\ack The authors are very grateful to Gregory Lupton for fruitful discussions and, in particular, for providing Example \ref{ex-lup}.

%%%%%%%%%%%%%%%%%%%%%%%%%%%%%%%%%% Section 2 %%%%%%%%%%%%%%%%%%%%%%%%%%%%%%%%%%%%%%

\section{Positive Quaternion K\"ahler Manifolds}\label{sec01}

Due to Berger's celebrated theorem the holonomy group $\Hol(M,g)$ of
a simply-connected, irreducible and non-symmetric Riemannian
manifold $(M,g)$ is one of $\SO(n)$, $\U(n)$, $\SU(n)$, $\Sp(n)$,
$\Sp(n)\Sp(1)$, $\G_2$ and $\Spin(7)$.

A connected oriented Riemannian manifold $(M^{4n},g)$ is called a
\emph{quaternion K\"ahler manifold} if
\begin{align*}
\Hol(M,g)\In \Sp(n)\Sp(1)=\Sp(n)\times \Sp(1)/\langle -\id, -1\rangle
\end{align*}
(In the case $n=1$ one additionally requires $M$ to be Einstein and
self-dual.) Quaternion K\"ahler manifolds are Einstein
(see \cite[Theorem 14.39, p.~403]{Bes08}). In particular, their scalar
curvature is constant.
\begin{defi}\label{PQKdef01}
A \emph{positive quaternion K\"ahler manifold} is a quaternion
K\"ahler manifold with complete metric and with positive scalar
curvature.
\end{defi}
For an elaborate depiction of the subject we recommend the survey
articles \cite{Sal82} and \cite{Sal99}. We shall content ourselves
with mentioning a few properties that will be of importance
throughout this article:

Foremost, we note that a positive quaternion K\"ahler manifold $M$
is  not necessarily K\"ahlerian, as the name might suggest.
Moreover, the manifold $M$ is compact and simply-connected
(see \cite[ p.~158]{Sal82} and \cite[6.6, p.~163]{Sal82}).

The only known examples of positive quaternion K\"ahler manifolds
are given by the so-called \emph{Wolf-spaces}, which are all
symmetric. They are the only possible homogeneous examples by a result of  Alekseevski.
They are given by the infinite series $\hh\pp^n$,
$\Gr_2(\cc^{n+2})$ and $\widetilde \Gr_4(\rr^{n+4})$ (the
Grassmanian of oriented real $4$-planes) and by the exceptional spaces
$\G_2/\SO(4)$, $\F_4/\Sp(3)\Sp(1)$, $\E_6/\SU(6)\Sp(1)$,
$\E_7/\Spin(12)\Sp(1)$, $\E_8/\E_7\Sp(1)$. Besides, it is known that
in each dimension there are only finitely many positive quaternion
K\"ahler manifolds (cf.~\cite{LS94}.0.1, p.~110). This endorses
the fundamental conjecture by LeBrun and Salamon (Conjecture \ref{conj01}) speculating that positive quaternion K\"ahler
manifolds are symmetric spaces.
% i.e.~Wolf-spaces indeed.

A confirmation of the conjecture has been achieved in dimensions
four (Hitchin) and eight (Poon--Salamon, LeBrun--Salamon). For a discussion of dimension $12$ see \cite{AD10} and \cite{HH10}.

%Thus one motivating question throughout this article will be: How close are positive quaternion K\"ahler manifolds to symmetric spaces? We shall
%give a partial answer from the viewpoint of Rational Homotopy Theory.

\vspace{5mm}

Remarkably, the theory of positive quaternion K\"ahler manifolds may be completely transcribed to an equivalent
theory in complex geometry. This is done via the \emph{twistor
space} $Z$ of the positive quaternion K\"ahler manifold $M$. This
Fano contact K\"ahler Einstein manifold may be constructed as follows:

Locally the principle  $\Sp(n)\Sp(1)$ structure bundle may be lifted
to its double covering with fiber $\Sp(n)\times \Sp(1)$. So, locally,
one may use the standard representation of $\Sp(1)$ on $\cc^2$ to
associate a vector bundle $H$. In general, $H$ does not exist
globally, but its complex projectivization $Z=\pp_\cc(H)$ does. In
particular, we obtain the \emph{twistor fibration}
\begin{align*}
\cc\pp^1\hto{} \pp_\cc(H) \to M
\end{align*}

Alternatively, the manifold $Z$ may be considered as the unit sphere
bundle $\s(E')$ associated to the $3$-dimensional subbundle $E'$ of
the vector bundle $\End(TM)$ generated locally by the almost complex
structures $I$, $J$, $K$ which behave like the corresponding unit
quaternions $i$, $j$ and $k$. That is, the twistor fibration is just
\begin{align*}
\s^2\hto{} \s(E') \to M
\end{align*}
(Comparing this bundle to its version above we need to remark that
clearly $\cc\pp^1\cong \s^2$.) The existence of this twistor bundle together with the fact that the total space is a compact K\"ahler manifold is basically the only property which we shall exploit in order to prove Theorem \ref{theoA}.

As an example one may observe that on $\hh\pp^n$ we have a global
lift of $\Sp(n)\Sp(1)$ and that the vector bundle associated to the
standard representation of $\Sp(1)$ on $\cc^2$ is just the
tautological bundle. Now complex projectivization of this bundle
yields the complex projective space $\cc\pp^{2n+1}$ and the twistor
fibration is just the canonical projection.

More generally, on Wolf spaces one obtains the following: The Wolf
space may be written as $G/K\Sp(1)$ (cf.~the table on \cite[p.~409]{Bes08}) and its corresponding twistor space is given as $G/K\U(1)$
with the twistor fibration being the canonical projection.

\vspace{5mm}

Using twistor theory a variety of remarkable results have been
obtained. Let us mention  just  the following ones:

%\begin{theo}[Strong rigidity]\label{PQKtheo01}
%Let $(M,g)$ be a Positive Quaternion K\"ahler
%Manifold. Then we have
%\begin{align*}
%\pi_2(M)=
%\begin{cases}
% 0 & \textrm{iff } M\cong \hh\pp^n\\
%\zz & \textrm{iff } M\cong \Gr_2(\cc^{n+2})\\
%\textrm{finite with } \zz_2\textrm{-torsion contained in } \pi_2(M)
%&\textrm{otherwise}
%\end{cases}
%\end{align*}
%\end{theo}
%\begin{prf}
%See theorems \cite{LS94}.0.2, p.~110, and \cite{Sal99}.5.5,
%p.~103.
%\end{prf}
\begin{theo}\label{PQKtheo02}
Odd-degree Betti numbers of $M$ vanish, i.e.~$b_{2i+1}=0$ for $i\geq 0$.
\end{theo}
\begin{prf}
See \cite[Theorem 6.6, p.~163]{Sal82}, where it is shown that the
Hodge decomposition of the twistor space is concentrated in terms
$H^{p,p}(Z,\rr)$.
\end{prf}
This implies that a rationally elliptic positive quaternion K\"ahler manifold is an $F_0$-space and that it is formal, in particular.

\begin{theo}[Strong rigidity]\label{PQKtheo01}
Let $(M,g)$ be a positive quaternion K\"ahler
manifold. Then we have
\begin{align*}
\pi_2(M)=
\begin{cases}
 0 & \textrm{iff } M\cong \hh\pp^n\\
\zz & \textrm{iff } M\cong \Gr_2(\cc^{n+2})\\
\textrm{finite with } \zz_2\textrm{-torsion contained in } \pi_2(M)
&\textrm{otherwise}
\end{cases}
\end{align*}
\end{theo}
\begin{prf}
See \cite[Theorem 0.2, p.~110]{LS94} and \cite[Theorem 5.5, p.~103]{Sal99}.
\end{prf}

%%%%%%%%%%%%%%%%%%%%%%%%%%%%%%%%%% Section 2 %%%%%%%%%%%%%%%%%%%%%%%%%%%%%%%%%%%%%%

\section{Rational Homotopy Theory}\label{sec02}

%%%%%%%%%%%%%%%%%%%%%%%%%%%%%%%%%%%%%%%%%%%%%%%%%%%%%%%%%%%%%%%%%%%%%%%%

\subsection{Formal spaces}

\begin{defi}\label{HistRHTdef03}
A commutative differential graded algebra $(A,\dif)$ (over a field
$\kk\Ni \qq$) is called \emph{formal}, if it is weakly equivalent to
the cohomology algebra $(H(A,\kk),0)$ (with trivial differential).

We call a path-connected topological space \emph{formal} if
% its rational homotopy type is a formal consequence of its rational
%cohomology algebra, i.e.~if
$(\APL(X),\dif)$ is formal.  In detail,
the space $X$ is formal if and only if there is a weak equivalence
$(\APL(X),\dif)\simeq (H^*(X),0)$, i.e.~a chain of
quasi-isomorphisms
\begin{align*}
(\APL(X),\dif)\xleftarrow{\simeq} \dots \xto{\simeq} \dots
\xleftarrow{\simeq} \dots \xto{\simeq} (H^*(X),0)
\end{align*}
\end{defi}
The algebras involved are algebras over the rationals. However, it turns out that the notion of formality does not depend on rational coefficients.
\begin{theo}\label{HistRHTtheo01}
Let $X$ have rational homology of finite type and let $\kk\Ni \qq$ be a field extension.

Then the algebra
$(\APL(X,\kk),\dif)$ is formal
if and only if $X$ is a formal space.
\end{theo}
\begin{prf}
See \cite[p.~156]{FHT01} and \cite[Theorem 12.1, p.~316]{FHT01}.
\end{prf}
Thus we need not worry about field extensions and it suffices to
consider rational coefficients only.
\begin{ex}\label{HistRHTex01}
The following spaces are formal:
\begin{itemize}
 \item
$H$-spaces (\cite[Example 12.3, p.~143]{FHT01}).
\item
symmetric spaces of compact type (\cite[Example 12.3,
p.~162]{FHT01}).
\item
$N$-symmetric spaces (see \cite[Main Theorem,
p.~40]{St02} for the precise statement, \cite{KT03}).
\item
compact K\"ahler manifolds (\cite[Main
Theorem, p.~270]{DGMS75}).
\end{itemize}
\end{ex}

\subsection{Formal maps}
Let us recall the following definition.
\begin{defin}{Formal maps}
Let $(A, \dif)$, $(A',\dif')$ be formal dga's and let $f\co (A,\dif)\to (A',\dif')$ be a morphism of dga's.
Let $\mu_A\co (M_A,\hat \dif)\to (A,\dif)$ and $\mu_A\co (M_{A'},\hat \dif')\to (A',\dif')$ be minimal Sullivan models.

Let $\hat f\co (M_A,\hat\dif)\to (M_{A'},\hat \dif')$ be a Sullivan representative, i.e.~an induced map of minimal models unique up to homotopy.
Then $f$ is called {\it formal} if there exist quasi-isomorphisms $m_A\co (M_A,\hat \dif)\to (H(A,\dif),0)$ and $m_{A'}\co (M_{A'},\hat\dif')\to (H(A',\dif'),0)$ which are the identity on cohomology and make the following diagram commute up to homotopy

 %
%%xybelow

\begin{equation}\label{formal-map}
\xymatrix{
(M_A,\hat\dif)\ar[r]^{\hat f}\ar[d]^{m_A}&(M_{A'},\hat\dif')\ar[d]^{m_{A'}}\\
(H(A,\dif),0)\ar[r]^{f^*}&(H(A',\dif'),0)
 }
\end{equation}
%%xyabove
%
where we identify $H(A,\dif)$ with  $H(M_A,\hat \dif)$ via $\mu_A^*$ and $H(A',\dif')$ with \linebreak[4] $H(M_{A'},\hat \dif')$ via $\mu_{A'}^*$ .

A map $f\co X\to Y$ between formal topological spaces is called formal if the induced map $f_\qq\co (A_{PL}(Y),\dif)\to (A_{PL}(X),\dif')$ is formal.
\end{defin}
\begin{rem2}\label{form-nonmin}
Let us remark that the minimality of $M_A$ and $M_{A'}$ is not necessary in the above definition and $f$ is formal if and only if one can find Sullivan models $M_A$ and $M_{A'}$ satisfying the above definition.
\end{rem2}

%%%%%%%%%%%%%%%%%%%%%%%%%%%%%%%%%%%%%%%%%%%%%%%%%%%%%%%%%%%%%%%%%%%%%%%%

\subsection{Absolute bigraded and filtered models}

\begin{theo}[Halperin--Stasheff bigraded model \cite{HaSt}]
Let $A$ be a finitely generated graded commutative algebra over $\qq$. We suppose that $A_0=\qq$. Then the cdga $(A, 0)$ admits a minimal model $\rho \co  (\Lambda V, \dif) \to (A, 0)$ where $V$ is equipped with a lower gradation $V = \oplus _{p\ge 0} V_p$ extended in a multiplicative way to $\Lambda V$ and where the following properties hold.
\begin{enumerate}
\item $\dif(V_p)\subset (\Lambda V)_{p-1}$. In particular, $\dif(V_0)=0$. Therefore the cohomology is a bigraded algebra $H(\Lambda V, \dif)=\oplus_{p\ge 0} H_p(\Lambda V,\dif)$.
\item  $\rho(V_p)=0$ for $p>0$.
\item  $H_p(\Lambda V,\dif)=0$ for $p>0$ and  $\rho^*\co H_0(\Lambda V,\dif)\to H(A,0)=A$ is an isomorphism
\end{enumerate}

The cdga $ (\Lambda V, \dif)$ is called a \emph{bigraded model} of the graded algebra $A$.
\end{theo}
 Bigraded models are unique in the natural sense.

As noted above the fact that $\dif$ is homogeneous of lower degree -1 means that $H_p(\Lambda V,\dif)\to H(\Lambda V,\dif)$ is injective for any $p$.
Also recall that the algebra of derivations of a dga $(C,\dif)$ denoted by $\Der(C)$ is naturally a differential graded Lie algebra with the differential $\dD$ given by the bracket with $\dif$:
\begin{align*}
\dD(\theta)=[\dif,\theta]=\dif\circ \theta- (-1)^q\theta\circ \dif
\end{align*}
where $\theta$ is a derivation of degree $q$. If $(C,\dif)=(\Lambda V, \dif)$ is a bigraded model of a dga $(A,0)$, the lower filtration on $\Lambda V$ induces a lower filtration on $\Der(\Lambda V,\dif)$. We denote by $\Der_{p}^q(\Lambda V,\dif)$ the set of derivations of $(\Lambda V,\dif)$ which increase the usual degree by $q$ and decrease the lower degree by $p$. With these notations we have that
\begin{align*}
\dD\co \Der_p^q(\Lambda V,\dif)\to \Der_{p+1}^{q+1}(\Lambda V,\dif)
\end{align*}
We shall denote the cohomology of this complex by $H_p^q(\Der(\Lambda V,\dif))$. As before, the $(-1)$-homogeneity of  $\dD$ with respect to the lower degree implies that the natural map $H_p^q(\Der(\Lambda V,\dif))\to H^{q}(\Der(\Lambda V,\dif))$ is injective for any $p,q$.

For cdga's with non-zero differentials Halperin and Stasheff (\cite{HaSt}) introduced so-called \emph{filtered models} exhibiting such algebras as deformations of their cohomology algebras in an appropriate sense.

\begin{theo}\label{filter}
Let $(A,\dif_A)$ be a connected cdga of finite type. Let \linebreak[4]$\rho\co (\Lambda V,\dif)\to (H(A,\dif_A),0)$ be a bigraded model of $A$.

Then there exists a differential $\operatorname{D}$ on $\Lambda V$ such that
\begin{enumerate}
\item $(\Lambda V,\operatorname{D})\rto{\pi} (A,\dif_A)$ is a Sullivan model of $(A,\dif_A)$.  (This model is not necessarily minimal.)
\item $\operatorname{D}-\dif$ decreases the lower degree by at least two:
\begin{align*}
\operatorname{D}-\dif\co V_p\to (\Lambda V)_{\le p-2}
\end{align*}
 That is, the differential $\operatorname{D}$ can be written as $\operatorname{D}=\dif+\dif_2+\dif_3+\ldots$ where $\dif_i$ is homogeneous of degree $-i$ in lower degree:
 \begin{align*}
 \dif_i\co V_p\to  (\Lambda V)_{p-i}
 \end{align*}
\end{enumerate}
\end{theo}

Filtered and bigraded models are useful for a number of reasons. Of  particular interest to us is the fact that they provide a good framework for distinguishing formal cdga's amongst all cdga's with a given cohomology ring.  As shown by Halperin and Stasheff, the deformation differentials $\dif_i$ define a series of obstructions $o_i$.

\begin{theo}[Obstructions to formality  \cite{HaSt}]\label{obstr-form}
Let $(A,\dif_A)$ be a connected cdga of finite type.  Let $\rho\co (\Lambda V,\dif)\to (H(A),0)$ be a bigraded model of $A$. Let $(\Lambda V,\operatorname{D})\rto{\pi} A$ with $\operatorname{D}=\dif+\dif_2+\dif_3+\ldots$ be  a filtered model of $A$ with respect to the chosen bigraded model.
Then we obtain:

If $\dif_j=0$ for $j<i$, then $\dif_i$ is a closed derivation in $\Der_i^1(\Lambda V,\dif)$ and its cohomology class in $H_i^1(\Der(\Lambda V,\dif))$  is denoted by $o_i$. If $o_i=0$, then there exists an automorphism of $\Lambda V$ as a cga (but not as  cdga) such that conjugating both the algebra and the differential by this automorphism yields a new filtered model $(\Lambda V,\operatorname{D}')$. This model has the important property that the decomposition of its differential
\begin{align*}
\operatorname{D}'=\dif'+\dif'_2+\dif'_3+\ldots
\end{align*}
as above satisfies that the $\dif'_j$ are identically zero for $j=2,3,\ldots, i$.

Repeating this process one obtains that
$A$ is formal if and only if the consecutive sequence of obstructions $o_i$ obtained in this fashion vanishes for all $i\ge 2$.
\end{theo}
This obstruction theory will be one of our main tools in proving Theorem \ref{theoA}.

\vspace{5mm}

Before we go on let us observe that the minimality of a bigraded model is not important in the above result. This fact is probably well-known; yet, since we do not have an explicit reference for this, we sketch a proof here.

Let $(A,0)$ be a cdga and let  $\phi\co (\Lambda V,\dif)\to (A,0)$ be a
%free ??????
Sullivan model of $A$ satisfying all the properties of a bigraded model except  possibly  minimality. Following Saneblidze \cite{Snblz} we shall refer to such a model as a {\it multiplicative resolution} of $A$. It is also a Tate-Jozefiak resolution of $A$ in the category of cgas (i.e.~forgetting the zero differential on $A$)---cf.~\cite{Joz72}.

In complete analogy to Theorem \ref{filter}, we may speak of a \emph{filtered model with respect to a multiplicative resolution}.

%Then theorem \ref{obstr-form} immediately implies

\begin{cor}\label{obstr-form-cor}
Let  $(A,\dif_A)$ be a connected cdga.  Let $\phi\co (\Lambda V,\dif)\to (H(A),0)$ be a multiplicative resolution of $A$. Let $\pi\co (\Lambda V,\operatorname{D})\rto{\simeq} (A,\dif_A)$ with $\operatorname{D}=\dif+\dif_2+\dif_3+\ldots$ be a filtered model of $(A,\dif_A)$ with respect to the multiplicative resolution $\phi$.

We derive the following results: If $\dif_j=0$ for $j<i$, then $\dif_i$ is a closed derivation in $\Der_{i}^1(\Lambda V,\dif)$ and its cohomology class in $H_{i}^1(\Der(\Lambda V,\dif))$  is denoted by $o_i$.

As a consequence, the cdga $(A,\dif_A)$ is formal if and only if $o_i=0$ for all $i\ge 2$.

\end{cor}
\begin{prf}
The ``if'' direction is proved in exactly the same way as in  \cite{HaSt} as it does not use minimality of the multiplicative resolution.
For the ``only if'' direction recall that Saneblidze constructed relative bigraded and filtered models for fibrations (see \cite{Saneblz91} or  \cite{Snblz}). While his construction is somewhat different from the ones in \cite{Tho81}  or \cite{V-P81}, which we discuss below, in the case the base equals a point it reduces to the filtered model of the fiber based on an arbitrary (not necessarily minimal) multiplicative resolution of the fiber. Moreover, Saneblidze also proved uniqueness of such filtered models (see \cite[Theorem 3.3]{Saneblz91} or \cite[Theorem A]{Snblz}) generalizing the uniqueness theorem for bigraded models of Halperin and Stasheff \cite[Theorem 4.4]{HaSt} where the same statement was proved for filtered models based on minimal multiplicative resolutions. This uniqueness implies the ``only if'' direction in exactly the same way as in  \cite{HaSt}.
\end{prf}

\subsection{Relative bigraded and filtered models}
We shall also need  relative versions of bigraded and filtered models developed by Vigu\'e--Poirrier  \cite{V-P81} (cf.~also \cite{Snblz, Saneblz91} and \cite{Tho82} ).

%We first define a bigraded model.
As in the absolute case, one begins by constructing a bigraded model in the relative category.
For this let $\phi \co (H,0)\to (H',0)$ be a morphism of two cdga's (with trivial differentials).
Let $\rho\co (\Lambda Z,\dif)\xto{\simeq} (H,0)$ be a bigraded model and let
\[
 (\Lambda Z,\dif)\hto{i} (\Lambda Z\otimes \Lambda X, \dif')\rto{p}(\Lambda X, \dif'')
\]
 be a minimal relative model for $\phi\circ \rho$. That is, $ (\Lambda Z,\dif)\hto{i} (\Lambda Z\otimes \Lambda X, \dif')$ is a minimal relative Sullivan algebra (see \cite[Section 14]{FHT} for definitions) such that the following diagram commutes.
\begin{align}\label{rel-bigr}
\xymatrix{
(H,0)\ar[r]^{\phi}&(H',0)\\
 (\Lambda Z,\dif)\ar[u]_\rho \ar@{-->}[ur]^{\phi\circ\rho} \ar[r]^<<<<{i} &(\Lambda Z\otimes \Lambda X, \dif')\ar[u]_{\rho'}\ar[r]^<<<<{p}&(\Lambda X, \dif'')
 }
\end{align}
Here both $\rho$ and $\rho'$ are quasi-isomorphisms. (Note that $(\Lambda Z\otimes \Lambda X, \dif')$ is a Sullivan algebra which is minimal as a relative algebra but might not be minimal as an absolute algebra.)

\begin{theo}[Relative bigraded model,  \cite{V-P81}]\label{rel-bigr-thm}
With the terminology from above one obtains the following: The  algebra $(Y,\dif'):=(\Lambda Z\otimes \Lambda X, \dif')$ can be chosen in such a way that it admits a lower grading satisfying the following properties

\begin{enumerate}[a.)]
\item $Z=\oplus_{i\ge 0} Z_i, X=\oplus_{i\ge 0} X_i, Y_i=X_i\oplus Z_i$.
\item The differential $\dif'$ is homogeneous of degree $-1$ with respect to the lower grading.
\item All the maps in diagram \eqref{rel-bigr} preserve the lower grading. For this $(H,0)$ and $(H',0)$ are understood to have trivial lower gradings, i.e.~$H_0=H$ and $H_i=0$ for $i>0$ respectively $H_0'=H'$ and $H_i'=0$ for $i>0$.
\item The morphism $\rho'\co (\Lambda Y,\dif')\to (H',0)$ is a multiplicative resolution. That is
\begin{enumerate}[(i)]
\item $\rho'\co \Lambda Y_0\to H'$ is onto.
\item $\rho'(Y_i)=0$ for $i>0$.
\item $H_i( \Lambda Y,\dif')=0$ for $i>0$ and $({\rho'})^*_0\co H_0(\Lambda Y,\dif')\to H'$ is an isomorphism.
\end{enumerate}
\end{enumerate}
\end{theo}

We shall refer to this relative Sullivan algebra
\begin{align*}
(\Lambda Z,\dif)\hto{i} (\Lambda Z\otimes \Lambda X, \dif')\rto{p}(\Lambda X, \dif'')
\end{align*}
as the \emph{(relative) bigraded model of $\phi$}.

Lastly, Vigu\'e--Poirrier (see  \cite{V-P81} ) proved the existence of a relative filtered model for morphisms between arbitrary  cdga's.

\begin{theo}\label{rel-filtered}
Let $\alpha\co (A,\dif)\to (A',\dif')$ be a morphism of cdga's with the property that $H^1(\alpha)$ is injective. Let
\begin{equation}\label{rel-fil-e1}
\xymatrix{
(H(A,\dif),0)\ar[r]^{\alpha^*}&(H(A',\dif'),0)\\
 (\Lambda Z,\dif)\ar[u]_\rho \ar@{-->}[ur]^{\alpha^*\circ\rho} \ar[r]^<<<<<<{i} &(\Lambda Z\otimes \Lambda X, \dif')\ar[u]_{\rho'}\ar[r]^<<<<{p}&(\Lambda X, \dif'')
 }
\end{equation}
be a bigraded model of $\alpha^*$ provided by Theorem \ref{rel-bigr-thm}.

Let $(\Lambda Z,\operatorname{D})\rto{\pi} (A,\dif)$ be a filtered model of $(A,\dif)$. Then $\alpha\circ\pi\co  (\Lambda Z,\operatorname{D})\to (A',\dif)$ admits a minimal relative model
\begin{align}\label{rel-fil-e2}
\xymatrix{
(A,\dif)\ar[r]^{\alpha}&(A',\dif')\\
 (\Lambda Z,\operatorname{D})\ar[u]_\pi \ar@{-->}[ur]^{\alpha\circ\pi} \ar[r]^<<<<{i} &(\Lambda Z\otimes \Lambda X, \operatorname{D}')\ar[u]_{\pi'}\ar[r]^<<<<{p}&(\Lambda X, \operatorname{D}'')
 }
\end{align}
where $\operatorname{D}'-\dif'$ decreases the lower grading by at least 2, i.e.
\begin{align}
\operatorname{D}'-\dif'\co (\Lambda Z\otimes \Lambda X)_p\to (\Lambda Z\otimes \Lambda X)_{\le p-2} \qquad\textrm{ for any } p\ge 0.
\end{align}

\end{theo}

We shall refer to the relative Sullivan  algebra
\begin{align*}
(\Lambda Z,\operatorname{D})\rto{i} (\Lambda Z\otimes \Lambda X, \operatorname{D}')\rto{p}(\Lambda X, \operatorname{D}'')
\end{align*}
as the \emph{(relative) filtered model of $\alpha$}.
\begin{rem}
In \cite{Tho81} Thomas considered a simplified version of a relative filtered model, which is obtained by pushing forward the above construction via the quasi-isomorphism $\pi\co (\Lambda Z, \operatorname{D})\to (A,\dif)$, i.e.~his filtered model can be obtained from the one defined above by taking $(A,\dif)\otimes_{(\Lambda  Z, \operatorname{D})} (\Lambda Z\otimes \Lambda X, \operatorname{D}')$.
\end{rem}

%%%%%%%%%%%%%%%%%%%%%%%%%%%%%%%%%%%%%%%%%%%%%%%%%%%%%%%%%%%%%%%%%%%%%%%%

\subsection{Bigraded and filtered models of TNCZ fibrations}

For a general filtered model of a map $\alpha\co A\to A'$ the fiber cdga $(\Lambda X, \operatorname{D}'')$ is not a filtered model of the fiber. However, it is one if the map $\alpha\co A\to A'$ is a model of a TNCZ fibration, which is the situation we are interested in in this article. More precisely, the following holds
\begin{theo}[\cite{V-P81} (cf.~\cite{Tho82})]\label{theo-tncz}
Let $F\hto{} E\rto{f} B$ be a Serre fibration of path-connected spaces where $B$ is simply-connected and $H^*(F)$ has finite type. Suppose this fibration is TNCZ.
Let $f_\qq\co M_B\to  M_E$ be the induced map of minimal models.  Then $f_\qq$ admits a relative filtered model
\begin{align*}
(\Lambda Z,\operatorname{D})\hto{i} (\Lambda Z\otimes \Lambda X, \operatorname{D}')\rto{p}(\Lambda X, \operatorname{D}'')
\end{align*}
such that  $(\Lambda X, \operatorname{D}'')$ is a filtered model of $F$.

In particular, $\operatorname{D}_1''=\dif''$  where $\dif''$ is the fiber differential in the bigraded model of  $f^*\co H^*(B)\to H^*(E)$. Here $\operatorname{D}''_1$ denotes the part of $\operatorname{D}''$ decreasing the lower degree by $1$.
%
%\[
%(\Lambda Z,\dif)\hto{i} (\Lambda Z\otimes \Lambda X, \dif')\rto{p}(\Lambda X, \dif'')
%\]
\end{theo}

%%%%%%%%%%%%%%%%%%%%%%%%%%%%%%%%%% Section 3 %%%%%%%%%%%%%%%%%%%%%%%%%%%%%%%%%%%%%%

\section{Proof of the main results}\label{sec03}

Before we  prove Theorem \ref{thm A'}, let us show by an example that the class of formal elliptic spaces satsifying the Halperin conjecture is strictly larger than the class of $F_0$-spaces satisfying Halperin's conjecture. In other words, there are formal elliptic spaces with vanishing Euler characteristic which do not possess non-trivial derivations of negative degree on their cohomology algebras.

However, we remark that a formal elliptic space is necessarily  two-stage due to \cite{FelHal82}. Moreover, as a consequence of finite dimensionality of cohomology, the filtration degree $1$ in this two-stage decomposition only consists of odd-degree elements. This follows from comparing the minimal model with the $E_0$-term of the \emph{odd spectral sequence} (see \cite[Section 32b, p.~438]{FHT01}), i.e.~the \emph{pure Sullivan algebra associated with the minimal model}. This associated algebra has finite cohomological dimension if and only if the minimal model does (see  \cite[Proposition 32.4, p.~438]{FHT01}). However, any even degree element in filtration degree one in the minimal model becomes an element in filtration degree zero in the associated pure algebra. Such an element freely generates an infinite sequence of non-vanishing cohomology classes of increasing degrees.

Thus, basically, the ``only difference'' of a formal elliptic space with an $F_0$-space, which admits a pure model, lies in the fact that in filtration degree zero odd degree elements may occur.

Nonetheless, this suffices to give
\begin{ex}\label{ex01}
We define a minimal Sullivan algebra $(\Lambda V,\dif)$ by
\begin{align*}
V=\langle a,b,c,d,u,v \rangle \qquad \textrm{with } &\deg a=\deg b=\deg c=\deg d=3,\\& \deg u=6, \deg v=11
\end{align*}
We set $\dif a=\dif b=\dif c =\dif d=\dif u=0$ and $\dif v=abcd+u^2$.

This algebra is easily seen to be the minimal model of a rationally non-trivial $\s^6$-bundle over $\s^3\times\s^3\times\s^3\times \s^3$. Therefore this algebra is elliptic and also formal, since the base is formal and the fiber is $F_0$. In order to compute the derivations of negative degree on its cohomology algebra, one notes that the cohomology algebra is generated by elements of degrees $3$ and $6$. Consequently, every non-trivial homogeneous such derivation has to have degree $-3$ or $-6$. However, it is impossible to specify a non-trivial derivation on $[v]$ and on the $[a],[b],[c],[d]$ (and extend it as a derivation to the whole of $H(\Lambda V,\dif)$) which would be compatible with the fact that $[abcd]=-[u^2]$ in cohomology.

As $(\Lambda V,\dif)$ has non-zero cohomology generators of odd degree, its spatial realisation is not an $F_0$-space.
\end{ex}
\newpage
%\begin{rem}
%It may well be that a formal elliptic space which satisfies Halperin's conjecture is necessarily $F_0$. Thus, it is possible %that Theorem \ref{thm A'} is actually  exactly equivalent to theorem \ref{theoA}.
%\end{rem}

\begin{proof}[\textsc{Proof of Theorem \ref{thm A'}}]
We are now ready to proceed with the proof of Theorem \ref{thm A'}.

Let $F\hto{} E\rto{f} B$ be a Serre fibration where $E,F,B$ are simply-connected and of finite type.  Suppose further that $H^*(F)$ is finite dimensional, $E$ and $F$ are  formal and that $F$ is a two-stage space satisfying Halperin's conjecture.

Let $f^*\co H^*(B)\to H^*(E)$ be the induced map on cohomology and let
%%xybelow
\begin{align}\label{rel-bigr1}
\xymatrix{
(H(B),0)\ar[r]^{f^*}&(H(E),0)\\
 (\Lambda Z,\dif)\ar[u]_\rho \ar@{-->}[ur]^{f^*\circ\rho} \ar[r]^(0.4){i} &(\Lambda Z\otimes \Lambda X, \dif')\ar[u]_{\rho'}\ar[r]^(0.55){p}&(\Lambda X, \dif'')
 }
\end{align}
%%xyabove
be its relative bigraded model.

Let $M_B$, $M_E$ be the minimal models of $B$ and $E$ and let $f_\qq\co M_B\to M_E$ be a Sullivan representative, a corresponding map of minimal models.

Let
\begin{align}\label{rel-filtr2}
\xymatrix{
M_B\ar[r]^{f_\qq}&M_E\\
 (\Lambda Z,\operatorname{D})\ar[u]_\pi \ar@{-->}[ur]^{f_\qq\circ\pi} \ar[r]^(0.4){i} &(\Lambda Z\otimes \Lambda X, \operatorname{D}')\ar[u]_{\pi'}\ar[r]^(0.55){p}&(\Lambda X, \operatorname{D}'')
 }
\end{align}
be the filtered model of $f_\qq$.

We decompose $\operatorname{D}=\dif+\dif_2+\dif_3+\ldots$ and $\operatorname{D}'=\dif'+\dif_2'+\dif_3'+\ldots$ as above. Note that, since
\begin{align*}
(\Lambda Z,\operatorname{D})\hto{i} (\Lambda Z\otimes \Lambda X, \operatorname{D}')\rto{p}(\Lambda X, \operatorname{D}'')
\end{align*}
is a relative Sullivan algebra, we have that
\begin{align}\label{rel-sul:eq}
\dif_i'|_{\Lambda Z\otimes 1}=\dif_i
\end{align}

Since $F$ satisfies the Halperin conjecture, the fibration $F\hto{} E\rto{f} B$  is TNCZ and therefore, by Theorem \ref{theo-tncz}, we obtain that $(\Lambda X, \dif'')$ is a bigraded model of $H(F)$ and $(\Lambda X, \operatorname{D}'')$ is a filtered model of $F$.
Since $F$ is formal and elliptic, by \cite{FelHal82} we know that $X_i=0$ for $i\ge 2$. Thus, for degree reasons, \begin{align}\label{pure:eq}
\operatorname{D}''-\dif''|_{1\otimes \Lambda X}=0
\end{align}
and, in particular, $\operatorname{D}''=\dif''$.

This is a very useful fact which simplifies the situation considerably as it means that the relative filtered model is completely determined by the relative bigraded model and the filtered model of the base.  In particular, it immediately implies the fact mentioned in the introduction that for a TNCZ fibration with elliptic and formal fiber the formality of the base space implies the formality of the total space.

Recall that $\dif_i\in\Der_i^1(\Lambda Z,\dif)$ and $\dif_i'\in \Der_i^1 (\Lambda Z\otimes \Lambda X, \dif')$. Since $\Der(-,-)$ is contravariant in the first argument and covariant in the second argument, there are no natural maps between $\Der_i^1(\Lambda Z)$ and $\Der_i^1 (\Lambda Z\otimes \Lambda X)$.
However, both of them naturally map to  $\Der_i^1(\Lambda Z,  \Lambda Z\otimes \Lambda X) $. Here we view  $(\Lambda Z\otimes \Lambda X,\dif')$ as a dga-module over $(\Lambda Z,\dif)$ via the dga-homomorphism $(\Lambda Z,\dif)\hto{} (\Lambda Z\otimes \Lambda X, \dif')$.
\begin{align}\label{rel-bigr1}
\xymatrix{
\Der_i^1(\Lambda Z,\Lambda Z )\ar[dr]^{j}&\\
& \Der_i^1(\Lambda Z,  \Lambda Z\otimes \Lambda X)\\
\Der_i^1 (\Lambda Z\otimes \Lambda X,\Lambda Z\otimes \Lambda X)\ar[ur]^{j'}&
 }
\end{align}
By \eqref{rel-sul:eq} and \eqref{pure:eq} we have that $j'(\dif'_i)=j(\dif_i)$ for every $i$.

The differential on $\Der_i^1(\Lambda Z,  \Lambda Z\otimes \Lambda X)$ is given by
\begin{align*}
\dD(\theta)=\dif'\circ \theta- (-1)^k\theta\circ \dif
\end{align*}
where $\theta\in \Der(\Lambda Z,  \Lambda Z\otimes \Lambda X) $ is a derivation of degree $k$. The morphisms $j$ respectively $j'$ commute with the differentials $\dif$ respectively $\dif'$ and $\dD$, since $(\Lambda Z\otimes \Lambda X,\dif')$ is a relative Sullivan algebra over $(\Lambda Z,\dif)$.

\vspace{5mm}

Let us assume that $\dif_j=0$ for $j<i$. Due to equations \eqref{rel-sul:eq} and \eqref{pure:eq} this equally  implies that $\dif_j'=0$ for $j<i$.

Then $\dif_i'$ is a closed derivation and, since $E$ is formal, by Corollary \ref{obstr-form-cor} it is also exact. Therefore, the cocycle $j'(\dif_i)\in  \Der_i^1(\Lambda Z,  \Lambda Z\otimes \Lambda X)$ is exact, too.

Our key observation is that the map $j$ is injective in cohomology.

\begin{lemma}\label{key-lem}
Under the above assumptions
\begin{align*}
j^*\co H^*(\Der (\Lambda Z,\Lambda Z ))\to H^*(\Der (\Lambda Z,  \Lambda Z\otimes \Lambda X) )
\end{align*}
is injective.
\end{lemma}
\begin{prf}
Consider the following filtration on $\Der (\Lambda Z,  \Lambda Z\otimes \Lambda X)$: Let
\begin{align*}
F^p(\Der (\Lambda Z,  \Lambda Z\otimes \Lambda X))
\end{align*}
be the set of derivations which increase the $\Lambda Z$-degree by at least $p$; note that $p$ can be negative.   This filtration is clearly invariant under the differential $\dD(\theta)=\dif'\circ \theta- (-1)^k\theta\circ \dif$ (with $\theta\in \Der(\Lambda Z,  \Lambda Z\otimes \Lambda X) $ of degree $k$).

Let us examine the spectral sequence arising from this filtration.  Note that the filtration is bi-infinite on every  $\Der_k(\Lambda Z,  \Lambda Z\otimes \Lambda X)$. Therefore it is not immediately clear why this spectral sequence converges. However, we compute
\begin{align*}
E_0&=\Der (\Lambda Z,  \Lambda Z\otimes \Lambda X )\cong \Der (\Lambda Z,  \Lambda Z)\otimes  \Lambda X
\intertext{and the differential $\dif_0$ on $E_0$ satisfies $\dif_0=\dif''$. It follows that}
%\begin{align*}
E_1&=\Der (\Lambda Z,  \Lambda Z)\otimes H(\Lambda V,\dif'')
\end{align*}
Since $H(F)\cong H(\bigwedge V,\dif'')$ is finite dimensional, this implies that all but finitely many rows of $E_1^{\ast,\ast}$ are zero. This easily implies that the spectral sequence converges (and does so after finitely many steps, i.e.~it collapses) to
\begin{align*}
H(\Der (\Lambda Z,  \Lambda Z\otimes \Lambda X))
\end{align*}
We readily compute that
\begin{align*}
E_2^{p,q}=H^p(\Der (\Lambda Z, \Lambda Z))\otimes H^q(\Lambda X, \dif'')
\end{align*}
and that $j^*$ is the edge homomorphism of this spectral sequence.

\vspace{5mm}

Recall that $H(\Lambda X, \dif'')\cong H(F)$.
By a standard  argument the multiplicative properties of the spectral sequence imply that non-trivial spectral sequence differentials $\dif_r$ (for $r\ge 2$) produce non-trivial derivations on $H(F)$ of negative degree.

More precisely, suppose that $r=2$ or that $\dif_t=0$ for $2\le t<r$ for $r\geq 3$. This identifies the $E_r$-term with the $E_2$-term.  Thus for $[\theta]\otimes [z]\in E_r^{p,q}$ with $[\theta]\in H(\Der  (\Lambda Z,  \Lambda Z))$ and $[z]\in H(\bigwedge V,\dif'')$  we write $\dif_r([\theta]\otimes [z])$ as
\begin{align*}
\dif_r([\theta]\otimes [z])=\Sigma_{i} [\theta_i]\otimes [z_i]
\end{align*}
where $\{[\theta_i]\}_{i\in I}$ is a homogeneous basis of $H(\Der  (\Lambda Z,  \Lambda Z))$ and the $[z_i]$ are the corresponding coefficients from $H(F)$.
% for $j\in J$ form a homogeneous basis of $H^*(F)$.
Write
\begin{align*}
\dif_r([\theta] \otimes [z]) =\Sigma_{i} [\theta_i]\otimes (\dif_r)_i([\theta]\otimes [z])
\end{align*}
%---the differentials $\dif_r$ for $r\geq 2$ are trivial on the base algebra.
where $(\dif_r)_i([\theta]\otimes [z])$ denotes the coefficient $[z_i]$ in $H(F)$ for $[\theta_i]$ of the term $\dif_r([\theta]\otimes [z])$.

Suppose $\dif_r$ is not identically zero so that
\begin{align*}
\dif_r([\theta]\otimes [z])\ne 0
\end{align*}
for some $[\theta]\in H(\Lambda X, \dif'')$ and $[z]\in H(F)$. Pick $i_0$ with  $ [\theta_{i_0}]\otimes [z_{i_0}] \ne 0$.
Let $\phi\co H(F)\to H(F)$ be given by
\begin{align*}
\phi([z])=(\dif_r)_{i_0}([\theta]\otimes [z])\neq 0
\end{align*}
Then the multiplicative properties of the spectral sequence imply that $\phi$ is a non-trivial derivation of $H(F)$ of negative degree.

\vspace{5mm}

Since by assumptions $F$ satisfies Halperin's conjecture, this means that the derivations spectral sequence degenerates on the $E_2$-term. This immediately implies that the edge homomorphism $j^*$ is injective.
\end{prf}

We can now easily finish the proof of Theorem \ref{theoA}. In the situation depicted in diagram \eqref{rel-bigr1} we observe the following for the deformation differentials $\dif_i$ from above: Since
\begin{align*}
j^*([\dif_i])={j'}^*([\dif_i'])={j'}^*(0)=0
\end{align*}
and since $j^*$ is injective by Lemma \ref{key-lem}, we have that
\begin{align*}
[\dif_i]=0\in H^1(\Der(\Lambda Z,\dif))
\end{align*}

Next, since $H_i^1(\Der(\Lambda Z,\dif))\hto{} H^1(\Der(\Lambda Z,\dif))$ is injective, we derive that $[\dif_i]=0$ as an element of $H_i^1(\Der(\Lambda Z,\dif))$.

By Theorem \ref{obstr-form} this  means that we can modify the filtered model $(\Lambda Z, \operatorname{D})$ so that $\dif_j=0$ for $j\le i$. This also modifies the relative filtered model so that using  \eqref{rel-sul:eq} and \eqref{pure:eq} we can also assume that   $\dif_j'=0$ for $j\le i$. Consequently, the differential $\dif'_{i+1}$ is closed---whence exact---and we may proceed by induction.

This results in the fact that all the obstructions $[\dif_i]$ vanish. Hence the space $B$ is formal due to Theorem \ref{obstr-form}.

\vspace{5mm}

Finally, in order to prove the last statement of the theorem, we observe that whenever  $E$ is formal, the above procedure yields a relative filtered model which coincides with a relative bigraded model.  This trivially implies that $f$ is formal by Remark \ref{form-nonmin}. The converse is also true by \cite[Proposition 2.3.4]{V-P81}, but we do not need it here.
\end{proof}
\begin{rem}
Examining the proof it is clear that it applies to any fiber $F$ which is of finite type, has finite dimensional cohomology, satisfies Halperin's conjecture and the bigraded model of which has height $2$. That is, $V_i(F)=0$ for any $i>1$.
In this case the $F$ is automatically formal, since all the deformation differentials $\dif_i$ (for $i\geq 2$) necessarily vanish for degree reasons.

By an inductive argument on the lower degree it is easy to see that a space of finite type with finite dimensional cohomology satisfies $\dim V_i<\infty$ for all $i$.
As a consequence, a formal space of finite type with finite dimensional cohomology and $V_i=0$ for $i>1$ is necessarily rationally elliptic.
\end{rem}
%%%%%%%%%%%%%%%%%%%%%%%%%%%%%%%%%% Section 6 %%%%%%%%%%%%%%%%%%%%%%%%%%%%%%%%%%%%%%

\section{Counterexamples}\label{sec04}

In view of Theorem \ref{thm A'} one may wonder if the assumptions on the fiber in that theorem can be weakened while retaining at least one direction of the theorem.
Since the  product of two spaces is formal if and only if so are the factors, it is clear that the assumption on the fiber being formal is needed.

However, this condition alone is  easily seen not to be sufficient and  some further restrictions on the fiber are clearly necessary in order to relate formality of the base and of the total space in a fibration.

In algebraic terms, the simplest examples of fibrations of simply-connected spaces in which either only both fiber and base are formal or only both fiber and total space are formal can be provided as follows: For the first case one may apply a ``degree-shift'' of generators in the respective minimal models of the bundle
\begin{align*}
\s^1\hto{} M^3 \to \s^1\times \s^1
\end{align*}
from \cite[p.~261]{DGMS75} where $M$ is the $3$-dimensional compact
Heisenberg manifold, which is considered the simplest non-formal compact
manifold. Now let the underlying vector space of a minimal algebra $(\Lambda V,\dif)$ be generated by elements $x,y,z$ with $\deg x=\deg y=n$, $\deg z=2n-1$, $n$ odd, $\dif x=\dif y=0$ and $\dif
z=xy$. Then it is a direct observation that it can be realized as the non-formal total space of an analogous fibration of simply-connected spaces.

For the second case construct the formal minimal algebra $(\Lambda V,\dif)$ with $V=\langle b,c,n\rangle$ and $\deg b=3, \deg c=4, \deg n=6$ and with $\dif b=\dif c=0$, $\dif n=bc$. We then realise it as the base space of a fibration with the formal total space provided by the relative minimal model $(\Lambda V\otimes \Lambda\langle z\rangle
,\dif)$ with $\dif z=c$, $\deg z=3$ and with fiber rationally an $\s^3$.

\vspace{5mm}

It is also easy to construct nice geometric examples of fiber bundles with elliptic fibers where both the base and the fiber are formal but the total space is not.
For example, let $E=\Sp(5)/\SU(5)$. This space is well-known to be non-formal. It fibers with fiber $\s^3$ over the biquotient $B=\Sp(1)\backslash \Sp(5)/\SU(5)$ where the action of $\Sp(1)$ on the left on $\Sp(5)$ comes from the embedding $\Sp(1)\xto{\rho} \Sp(5)$ with $\rho(g)=\diag(g,1,\ldots, 1)$. (It is easy to see that the resulting action of $\Sp(1)\times \SU(5)$ on $\Sp(5)$ is free). Notice that the biquotient $B$ is $F_0$ and hence formal.
Thus we have a fibration $\s^3\hookrightarrow E\to B$ with formal base but non-formal total space.

If one does not insist on a simply-connected fiber, the fibration
\begin{align*}
\s^1\hto{}\Sp(n)/\SU(n) \to \Sp(n)/\U(n)
\end{align*}
for $n\geq 5$ is an even simpler example of that kind.

Note that there are many non-formal homogeneous spaces $G/H$ (such as aforementioned $\Sp(5)/\SU(5)$ or $\SU(6)/(\SU(3)\times \SU(3))$). Any such space fits into a fibration $H\hookrightarrow G\to G/H$ where both  the fiber and the total space are formal (and elliptic), but the base space is not.

\vspace{5mm}

As was mentioned earlier, Thomas constructed an example of a fibration $\s^3\vee \s^3\hookrightarrow E\to \s^3\times \s^5$ with $H(E)\cong H(B)\otimes H(F)$ as an algebra, but where $E$ is not formal (see \cite[Example II.13]{Tho82}). Note that here both $B$ and $F$ are obviously formal and the fibration is TNCZ.

The following example is due to Greg Lupton.
\begin{ex}\label{ex-lup}
We shall produce a fibration with formal, yet hyperbolic fiber $F$, which satisfies the Halperin conjecture and with a formal base space $B=\s^3$. However, the total space of this fibration is not formal.

Let $F=\s^2\vee \s^2\vee \s^2$ and $B=\s^3$. Then $F$ is formal and its bigraded model $(\Lambda V,\dif)$ is as follows:
\begin{align*}
V_0&=\langle a,b,c\rangle \\
V_1&=\langle\alpha,\beta,\gamma,\delta,\varepsilon,\phi\rangle\\
V_2&=\langle w,\ldots\rangle\\
\vdots
\end{align*}
where $\deg a=\deg b=\deg c=2$, $\dif a=\dif b=\dif c=0$, $\deg \alpha=\deg \beta=\deg \gamma=\deg \delta=\deg \varepsilon=\deg \phi=3$, $\dif\alpha=a^2, \dif\beta=ab, \dif\gamma=c^2, \dif\delta=b^2, \dif\varepsilon=ac, \dif\phi=bc$, $\deg w=4$, $\dif w=\alpha b-a\beta$, etc.

We will construct a closed derivation $\theta$ of $\Lambda V$ decreasing the lower and the upper degrees by $2$.
Define $\theta$ to be $0$ on $V_0$ and $V_1$, $\theta(w)=c$ and $\theta=0$ on the rest of the generators of $V_2$.
It is trivial to check that $[\dif, \theta]=0$  on $(\Lambda V)_{\leq 2}$.

We claim that $\theta$ can be extended to $\Lambda V$ to be a closed derivation of bi-degree $(2,-2)$.

We proceed by induction on lower degree.
Assume $i\ge 3$ and we have constructed $\theta$ on $(\Lambda V)_{\leq i-1}$  so that $[\dif, \theta]=0$ on $(\Lambda V)_{\leq i-1}$.
We claim that we can extend $\theta$ to  $(\Lambda V)_{\leq i}$  so that $[\theta, \dif]=0$ on $(\Lambda V)_{\leq i}$.

In order to prove this claim, we pick a basis of $V_i$ and let $v\in V_i$ be an element of this basis.
We have that $\dif v\in (\Lambda V)_{i-1}$ and due to the induction hypothesis
\begin{align*}
0=[\dif,\theta](\dif v)=\dif(\theta(\dif v))
\end{align*}
so that $\theta(\dif v)\in (\Lambda V)_{i-3}$ is closed. We then obtain two cases we need to consider separately.

\noindent\case{1} First suppose $i=3$.

Note that $\deg \dif v>4$ and hence $\deg \theta(\dif v)>2$. Since $H_0(\Lambda V,d)\cong H^*( \s^2\vee \s^2\vee \s^2)$ is zero in usual degree $>2$ we have that $\theta (\dif v)$ is exact. Therefore there is an $x\in (\Lambda V)_1$ with $\dif x=\theta(\dif v)$. Set $\theta(v)=x$. Then $[ \dif, \theta ](v)=0$.

\noindent\case{2}  Now suppose $i\ge 4$.

Then we have that $i-3\ge 1$. Since $H_{>0}(\Lambda V,\dif )=0$, this again implies that $\theta(\dif v)$ is exact. As before  we choose $x\in (\Lambda V)_{i-2}$ with $\dif x=\theta(\dif v)$ and set $\theta(v)=x$. Then again $[ \dif, \theta ]((v)=0$.

Thus, in any case we can extend the derivation to a closed derivation on $(\Lambda V)_{\leq i}$.

\vspace{5mm}

Next notice that $[\theta]\ne 0$ in $H_2^{-2}(\Der(\Lambda V,\dif))$. Indeed, suppose $\theta=[\dif, \mu]$ for some $\mu\in \Der_{1}^{-3}(\Lambda V,\dif)$. For degree reasons we must have that $\mu=0$ on $V_0$. Then
\begin{align*}
c=\theta(w)=[\dif,\mu](w)=\dif (\mu(w))-\mu( \dif w)=\dif (\mu(w))-\mu(\alpha b-a\beta)
\end{align*}
is in the ideal $I(a,b)$ (generated by $a$ and $b$) up to a coboundary. However,  $c\notin I(a,b)$ up to coboundary and hence $\theta$ is not exact.

Recall that we have an isomorphism $H^{-2}(\Der(\Lambda V,\dif))\cong \pi_2(\Aut_1 F)\otimes \qq$ constructed as follows.
For any formal space $W$  and any "nice" space $X$ (e.g.  simply connected and homotopy equivalent to a CW complex) we have a bijection
\begin{align*}
[W,X_\qq]\cong[M_X,M_W]\cong[M_X, H^*(W,\qq)]
\end{align*}
where $M_X$ and $M_W$ are minimal models of $X$ and $W$ respectively. Applying this to $W=F_\qq\times \s^2$ and $X=F_\qq$ we get
\begin{align*}
[F_\qq\times \s^2, F_\qq]&\cong [(\Lambda V,\dif), H^*(F)\otimes H(\s^2)]
\\&\cong  [(\Lambda V,\dif), (\Lambda V,\dif)\otimes H^*(\s^2)]
\\&\cong  [(\Lambda V,\dif), (\Lambda V,\dif)\otimes \Lambda \langle u\rangle/ (u^2)]
\end{align*}
where we identified $H^*(\s^2)$ with $\Lambda \langle u\rangle/(u^2)$ with $\deg u=2$ and  zero differential in the last equality.

Consider the map $h\co (\Lambda,\dif)\to (\Lambda,\dif)\otimes \Lambda (u)/(u^2)$ given by $h(x)=x\otimes 1+\theta(x)\otimes u$. The fact that $\theta$ is a closed derivation immediately implies that this is a dga homomorphism. Let $\tilde h\co F_\qq\times \s^2\to F_\qq$ be the corresponding element of $[F_\qq\times \s^2, F_\qq]$. By the adjunction formula it defines an element
\begin{align*}
\bar h\in \pi_3(\B\Aut_1(F))\otimes \qq \cong \pi_2(\Aut_1(F))\otimes \qq\cong \pi_2(\Aut_1(F_\qq))
\end{align*}

Let $F\hto{} E\xto{f} \s^3$ be the pullback  via $\bar h$ of the universal fibration
\begin{align*}
F \hto{}  \B\Aut_1^\bullet(F)\to \B\Aut_1(F)
\end{align*}
where $\Aut_1^\bullet(F)$ is the monoid of self-homotopy equivalences of $F$ homotopic to the identity relative to the base point.
By construction, the minimal model of $f$ is given by
\begin{align*}
(\Lambda \langle v\rangle,0)\hto{} (\Lambda\langle v\rangle\otimes \Lambda V,\operatorname{D})\to (\Lambda V, \dif)
\end{align*}
where $\deg v=3$ and the lower degree of $v$ is $0$.  Also, $\operatorname{D}(v)=0, \operatorname{D}(x)=\dif x+v\theta(x)$ for any $x\in \Lambda V$.
We can view $(\Lambda\langle v\rangle \otimes \Lambda V,\operatorname{D})$ as an absolute filtered model of $E$. Observe that $\operatorname{D}(w)=\dif w+vc$. Therefore, the $\dif_2$-part of $\operatorname{D}$ is not zero, since the lower degree of $w$ is 2 and the lower degree of $vc$ is $0$.

By the same argument as before it is easy to see that $\dif_2$ cannot be exact. Therefore the obstruction class $o_2\ne 0$ and hence $E$ is not formal.
On the other hand, the fiber and the base are  formal and the cohomology algebra of the fiber obviously has no non-trivial negative degree derivations; so, in particular, the fibration  $F\hto{} E\xto{f} \s^3$ is TNCZ.

%we finish with the following
\end{ex}
Motivated by Example \ref{ex-lup} it remains natural to ask for a fibration $F\hto{} E \xto{f} B$ with formal fiber $F$ and formal total space $E$ and with $F$ satisfying Halperin's conjecture, whilst $B$ is not formal.

%\begin{ques}
%Does there exist a fibration $F\hto{} E\xto{f} B$ where $F$ and $E$ are formal, $F$ satisfies Halperin's conjecture but $B$ %is not formal?
%\end{ques}

%Note however, that  we do not know of such examples where the fibration is  also TNCZ.

%%%%%%%%%%%%%%%%%%%%%%%%%%%%%%%%%% Bibliography %%%%%%%%%%%%%%%%%%%%%%%%%%%%%%%%%%%

%\bibliography{master-new}
%\bibliographystyle{abbrv}
\bibliographystyle{alpha}

\pagebreak \

\vfill

\begin{center}
\noindent
\begin{minipage}{\linewidth}
\small \noindent \textsc
{Manuel Amann} \\
\textsc{Department of Mathematics}\\
\textsc{University of Toronto}\\
%\textsc{Earth Sciences 2146}\\
\textsc{Toronto, Ontario}\\
\textsc{M5S 2E4} \\
\textsc{Canada}\\
[1ex]
\textsf{mamann@uni-muenster.de}\\
\textsf{http://individual.utoronto.ca/mamann/}
\end{minipage}
\end{center}

\vspace{10mm}

\begin{center}
\noindent
\begin{minipage}{\linewidth}
\small \noindent \textsc
{Vitali Kapovitch} \\
\textsc{Department of Mathematics}\\
\textsc{University of Toronto}\\
\textsc{Toronto, Ontario}\\
\textsc{M5S 2E4} \\
\textsc{Canada}\\
[1ex]
\textsf{vtk@math.toronto.edu}\\
\textsf{http://www.math.toronto.edu/vtk/}
\end{minipage}
\end{center}

\end{document}